\documentclass[a4paper,reqno,11pt]{amsart}
\usepackage[T1]{fontenc}
\usepackage[utf8]{inputenc}
\usepackage[foot]{amsaddr}    
\usepackage[english]{babel}
\usepackage{csquotes}         
\usepackage[babel]{microtype} 
\usepackage{dsfont}
\usepackage{amsmath}
\usepackage{amssymb}
\usepackage{tikz}
\usetikzlibrary{calc}
\usepackage{mathtools}
\usepackage{comment}

\usepackage{amsthm}
\usepackage{amstext}          
\usepackage{mleftright}       
\usepackage{fullpage}
\usepackage{xcolor}
\usepackage{caption}          

\usepackage{stickstootext}
\usepackage[stickstoo,varbb]{newtxmath} 
\makeatletter
\DeclareFontFamily{U}{ntxmia}{\skewchar \font =127}
 \DeclareFontShape{U}{ntxmia}{m}{it}{
                        <-> \ntxmath@scaled ntxmia
                      }{}    
                      \DeclareFontShape{U}{ntxmia}{b}{it}{
                        <-> \ntxmath@scaled ntxbmia
                      }{}
\makeatother

\usepackage{enumitem}         

\usepackage[numbers,sort]{natbib}

\makeatletter
\def\NAT@spacechar{~}
\makeatother

\usepackage{hyperref}              
\usepackage[capitalise]{cleveref}  
\hypersetup{colorlinks=true,
   citecolor=blue,
   filecolor=blue,
   linkcolor=blue,
   urlcolor=blue
}
\usepackage{url}

\allowdisplaybreaks

\makeatletter
\if@cref@capitalise
\crefname{figure}{figure}{figures}
\crefname{claim}{Claim}{Claims}
\crefname{conjecture}{Conjecture}{Conjectures}
\else
\crefname{figure}{Figure}{Figures}
\crefname{claim}{claim}{claims}
\crefname{conjecture}{conjecture}{conjectures}
\fi
\makeatother
\Crefname{figure}{Figure}{Figures}
\Crefname{claim}{Claim}{Claims}
\crefname{conjecture}{Conjecture}{Conjectures}

\usepackage{algorithm}
\usepackage{algpseudocode}
\newtheorem{theorem}{Theorem}
 
\newtheorem*{lemma*}{Lemma}

\theoremstyle{remark}

\newtheorem{example}[theorem]{Example}

\newtheorem{claim}[theorem]{Claim}

\theoremstyle{definition}

\newtheorem{conjecture}[theorem]{Conjecture}

\newcommand{\bw}{\operatorname{bw}}
\newcommand{\eps}{\epsilon}

\newcommand*\nume{\ensuremath{\mathrm{e}}}

\newcommand{\oldqed}{}
\def\endofClaim{\hfill\scalebox{.6}{$\blacksquare$}}

\newenvironment{claimproof}[1][Proof]{
	\renewcommand{\oldqed}{\qedsymbol}
	\renewcommand{\qedsymbol}{\endofClaim}
	\begin{proof}[#1]
	}{
	\end{proof}
	\renewcommand{\qedsymbol}{\oldqed}
}

\newcommand{\COMMENT}[1]{}

\title{Dirac's theorem for graphs of bounded bandwidth}

\author{Alberto Espuny Díaz}
\address[Espuny Díaz]{Institut f\"ur Informatik, Universit\"at Heidelberg, 69120 Heidelberg, Germany.}\email{espuny-diaz@informatik.uni-heidelberg.de}

\author{Pranshu Gupta}
\address[Gupta]{Faculty of Computer Science and Mathematics, University of Passau, Innstraße 33, Passau, Germany.}\email{pranshu.gupta@uni-passau.de}

\author{Domenico Mergoni Cecchelli}
\address[Mergoni Cecchelli]{Department of Mathematics, London School of Economics, Houghton Street, London WC2A 2AE, UK.}\email{d.mergoni@lse.ac.uk}

\author{Olaf Parczyk}
\address[Parczyk]{Zuse Institute Berlin, Department AIS2T, Takustraße 7,
14195 Berlin, Germany.}\email{parczyk@zib.de}

\author{Amedeo Sgueglia}
\address[Sgueglia]{Department of Mathematics,
University College London,
25 Gordon Street,
London, WC1H 0AY, UK.}\email{a.sgueglia@ucl.ac.uk}

\thanks{A.~Espuny Díaz has been funded by the Deutsche Forschungsgemeinschaft (DFG, German Research Foundation) through project no.\ 513704762.
O.~Parcyzk has been funded by the DFG under Germany's Excellence Strategy - The Berlin Mathematics Research Center MATH+ (EXC-2046/1, project ID: 390685689).
A.~Sgueglia has been funded by the Royal Society.
}
\usepackage{pgfplots}
\pgfplotsset{compat=1.18}

\begin{document}

\begin{abstract}
We provide an optimal sufficient condition, relating minimum degree and bandwidth, for a graph to contain a spanning subdivision of the complete bipartite graph $K_{2,\ell}$.
This includes the containment of Hamilton paths and cycles, and has applications in the random geometric graph model.
Our proof provides a greedy algorithm for constructing such structures.
\end{abstract}

\maketitle

\section{Introduction}

A classical theorem of \citet{Dirac52} asserts that every graph on $n\geq3$ vertices with minimum degree at least $n/2$ contains a Hamilton cycle.
This began a long line of research into sufficient minimum-degree conditions for (hyper/di)graphs to contain different (almost) spanning structures; see, e.g., the surveys of \citet{KO09Sur,KO12,KO14}, of \citet{RR10}, and of \citet{Zhao16}.
Moreover, the search for sufficient conditions for Hamiltonicity has been a driving force in graph theory; we refer the reader to the surveys of \citet{Gould91,Gould03,Gould14}.

In recent work on local resilience of random geometric graphs, \citet{ELW24} posed the following interesting conjecture, concerning a sufficient minimum-degree condition for subgraphs of $C_n^k$, the $k$-th power of a cycle $C_n$ on $n$ vertices, to contain a Hamilton cycle.\footnote{Here, the $k$-th power~$G^k$ of a graph~$G$ is obtained from~$G$ by adding an edge between any pair of vertices which are at graph distance at most~$k$ from each other.}

\begin{conjecture}[{\cite[Conjecture~1.13]{ELW24}}]\label{conj:cycle}
For all integers\/ $n\geq3$ and\/ $k\in [1,n/2]$, every graph\/ $G\subseteq C_n^k$ with\/ $\delta(G)\geq k+1$ is Hamiltonian.
\end{conjecture}

Motivated by \cref{conj:cycle}, we investigate the analogous question when the host graph is $P_n^k$, the $k$-th power of a path $P_n$ on $n$ vertices.
For an $n$-vertex graph $G$, the condition that $G \subseteq P_n^k$ is equivalent to $G$ having \emph{bandwidth} at most $k$, that is, that there exists a labelling $v_1,\ldots,v_n$ of $V(G)$ such that for every edge $v_iv_j\in E(G)$ we have that $|i-j|\leq k$.
The bandwidth of $G$, denoted by $\bw(G)$, is the smallest integer $k$ such that $G$ has bandwidth at most $k$.
The problem of determining the bandwidth of graphs, rooted in applications in computer science, has prompted a lot of research over the years (see, e.g., the surveys of \citet{CCDG82} and \citet{LW99}). Remarkably, graphs of bounded bandwidth are amenable for embedding problems in extremal graph theory~\cite{BST09,CKKO19,ABET20}.

Through the equivalence outlined above, our goal is to understand the interplay between minimum degree and bandwidth as sufficient conditions for subgraph containment.
However, the usual notion of minimum degree is not suitable in our context: while every vertex of $C_n^k$ has degree $2k$, there are vertices in $P_n^k$ of lower degree.
Thus, we have to adjust the notion of minimum degree accordingly.

Given a host graph $H$ and a spanning subgraph $G\subseteq H$, we say that $G$ has \emph{effective} minimum degree at least $\ell$ (with respect to $H$) if, for every $v\in V(H)$, we have that $\deg_{G}(v)\geq \min \{ \ell, \deg_H(v)\}$.
We let $\delta_\nume^H(G)$ denote the maximum $\ell$ such that $G$ has effective minimum degree at least $\ell$ with respect to~$H$.
In our problem, the actual embedding of a graph $G$ with $\bw(G)\leq k$ into $P_n^k$ is not relevant.
Thus, for the sake of conciseness, for an $n$-vertex graph $G$, we write $\delta_\nume^k (G)\geq \ell$ to indicate that there is an embedding of $G$ into $P_n^k$ such that $\delta_\nume^{P_n^k}(G)\geq \ell$.
Note that, in particular, this implies that $\bw(G) \le k$.

Our first contribution is an analogue of Dirac's theorem for graphs of bounded bandwidth.

\begin{theorem}\label{thm:cycle}
Let\/ $k$ and $n$ be integers with\/ $k\geq2$ and\/ $n\geq4$.
Any $n$-vertex graph $G$ with $\delta^k_\nume(G)\geq k+2$ contains a Hamilton cycle.
\end{theorem}

Note that, if $n$ is much larger than $k$, then `most' vertices in $P_n^k$ have degree $2k$.
In analogy to Dirac's theorem, \cref{thm:cycle} shows that, even after reducing the degree of most vertices by almost half, we can still guarantee a Hamilton cycle.
Observe that, for $k<n\leq2k$, \cref{thm:cycle} holds by Dirac's theorem, and for $n\in\{2k+1,\ldots,2k+4\}$ it holds, e.g., by the sufficient degree-sequence condition for Hamiltonicity of \citet{Posa62}.
Moreover, while it is stated for $n \ge 4$, it trivially also holds when $n=3$.
Finally, we remark that, for any $k \ge 2$, the condition on $\delta_\nume^k(G)$
in \cref{thm:cycle} is best possible for all $n\geq2k+3$ (see the construction in \cref{example:construction}).

\Cref{thm:cycle} is a special case of a more general result for subdivisions.
A \emph{subdivision} of a graph $H$ is obtained by replacing each edge of $H$ by a path of some positive length, all such paths being internally disjoint.
The problem of determining sufficient minimum-degree conditions for the containment of spanning subdivisions of different graphs has recently been considered by \citet{PS24} and \citet{Lee23}.
To state our result in the setting of graphs of bounded bandwidth, we need to consider rooted subdivisions.
Given a graph $G$ and two distinct vertices $u,v\in V(G)$, we say that $G$ contains a subdivision of $K_{2,\ell}$ \emph{rooted} at $u$ and $v$ if it contains such a subdivision where the maximal independent set of size $2$ in $K_{2,\ell}$ is embedded into $\{u,v\}$.

\begin{theorem}\label{thm:general}
    Let\/ $n$, $k$ and $\ell$ be integers with\/ $k\geq\ell\geq1$ and\/ $n\geq\ell+2$.
    Any $n$-vertex graph $G$ with $\delta^k_\nume(G)\geq k+\ell$ contains a spanning subdivision of\/ $K_{2,\ell}$ rooted at its two vertices of degree\/ $k$.\COMMENT{In fact, the subdivision may be rooted at any of the first $\ell+1$ vertices and any of the last $\ell+1$ vertices.}
\end{theorem}

Observe that the case $\ell=1$ corresponds to a Hamilton path and the case $\ell=2$ corresponds to \cref{thm:cycle}.
The next construction shows that the effective minimum degree condition in \cref{thm:general} (and in \cref {thm:cycle}) cannot be improved if $n$ is sufficiently large.

\begin{example}
\label{example:construction}
	Let $k\geq\ell\geq1$ and $n\geq2k+\ell+1$, and let $v_1,\ldots,v_n$ be $n$ distinct vertices.
	Let $G_1$ be the $k$-th power of the path $v_1v_2\ldots v_{k+\ell}$ and $G_2$ be the $k$-th power of the path $v_{k+2}v_{k+3}\ldots v_n$.
    Let~$G$ be the union of~$G_1$ and~$G_2$.
    Then, $G$ is an $n$-vertex graph and $v_1,\ldots,v_n$ is a labelling of $V(G)$ which witnesses that $\bw(G)\leq k$.
    Since both paths above have at least $k+\ell$ vertices, it follows from the construction that $\delta^k_\nume(G)=k+\ell-1$.
    Moreover, clearly there are only two vertices of degree $k$, which are $v_1$ and $v_n$, and removing the $\ell-1$ vertices $v_{k+2},\ldots,v_{k+\ell}$ disconnects $v_1$ and $v_n$.
    Thus, by Menger's theorem, $G$ cannot contain a subdivision of $K_{2,\ell}$ rooted at $v_1$ and $v_n$.
    Note that, when $\ell\in\{1,2\}$, $G$~cannot contain any spanning subdivision of $K_{2,\ell}$ at all (while, for other values of $\ell$, $G$ could contain a spanning subdivision of $K_{2,\ell}$ rooted at vertices different from $v_1, v_n$).\COMMENT{Say that $n$ is large. The subdivisions we obtain satisfy that every path has length $\Omega(n/k)$. The subdivisions that one could still potentially obtain with smaller minimum degree conditions (so, not rooted at the endpoints) would have paths of length $O(1)$.}
\end{example}    

Dirac’s theorem has been strengthened in several ways.
For example, \citet{Posa62} proved a condition on the degree sequence of a graph forcing Hamiltonicity.
In analogy to that strengthening, we consider conditions on the degree sequence (as opposed to effective minimum degree) in graphs of bounded bandwidth and obtain a stronger version of \cref{thm:general}.

\begin{theorem}\label{thm:deg-seq}
    Let\/ $n$, $k$ and $\ell$ be integers, with\/ $k\geq\ell\geq1$ and\/ $n\geq\ell+2$.
    Let\/ $G$ be an\/ $n$-vertex graph with\/ $\bw(G)\leq k$, and let\/ $v_1,\ldots,v_n$ be a labelling of\/ $V(G)$ witnessing this fact.
    Suppose that, for each\/ $i\in[n]$,\/ $\deg_G(v_i)\geq\min\{\ell+i-1,k+\ell,k+n-i\}$.
    Then,\/ $G$ contains a spanning subdivision of\/ $K_{2,\ell}$ rooted at\/ $v_1$ and\/ $v_n$.\COMMENT{In fact, the subdivision may be rooted at $v_1$ and any of the last $\ell+1$ vertices.}
\end{theorem}

Observe that the only difference between \cref{thm:general,thm:deg-seq} lies in the assumption on the degree of the first (at most) $k$ vertices (in the labelling witnessing that $\bw(G) \le k$): if $n\geq k+\ell-1$, in \cref{thm:general}, we require that $\deg(v_i) = k+i-1$ for $i \in [\ell]$ and that $\deg(v_i) \geq k+\ell$ for $i \in [k]\setminus[\ell]$, while in \cref{thm:deg-seq} we only require that $\deg(v_i) \ge \ell+i-1$ for $i \in [k]$.
The degree-sequence condition in \cref{thm:deg-seq} cannot be improved for any $i\in[k-\ell]$.
We refer to \Cref{fig:compare} for an easier comparison between the degree conditions forcing the existence of a Hamilton cycle in \cref{thm:general,thm:deg-seq}.

We remark that the proof of \cref{thm:deg-seq} (and thus of \cref{thm:general,thm:cycle}) is constructive and provides an efficient greedy algorithm to construct such a spanning subgraph, provided a labelling of $G$ witnessing its bandwidth is known.

\subsection{Local resilience and random geometric graphs}

Our concept of effective minimum degree is closely related to the problem of \emph{local resilience}, which examines the degree conditions under which a subgraph of a given graph maintains one of its properties.
This research direction, systematically initiated by \citet{SV08}, typically considers the problem where each vertex retains a proportion of its original degree.
The effective minimum degree is a variant that imposes an absolute lower bound on the degrees of all vertices, except those with initially lower degrees.
In this sense, our result addresses the local resilience of powers of paths with respect to Hamiltonicity and the containment of spanning subdivisions of $K_{2,\ell}$.
For a summary of all results on local resilience of random graphs with respect to Hamiltonicity, see the annotated bibliography of \citet{FriezeBiblio}.

Additionally, our results have applications on the local resilience of \emph{$1$-dimensional random geometric graphs}, which are defined as follows.
Given an integer $n$ and a real number $r\in[0,1]$, we let $G(n,r)$ denote a graph sampled by placing $n$ points independently and uniformly at random in $[0,1]$ and joining any pair of them by an edge if their distance is at most $r$.

By a simple concentration argument (see~\cite[Remark~1.15]{ELW24}), it follows that, for every $\eps>0$, there exists a constant $C>0$ such that, if $r\geq C\log n/n$, then a.a.s.\ $P_n^{(1-\eps/3)nr} \subseteq G(n,r) \subseteq P_n^{(1+\eps/3)nr}$.
Then, \cref{thm:cycle} immediately implies that a.a.s.\ every graph $H\subseteq G(n,r)$ with $\delta_\nume^{G(n,r)}(H)\geq(1+\eps)nr$ contains a Hamilton cycle (and, in fact, by \cref{thm:general}, a spanning subdivision of $K_{2,\ell}$ for any fixed~$\ell$).
This makes some progress towards the $1$-dimensional case of a conjecture of \citet[Conjecture~1.3]{ELW24}, who conjectured that the same result holds for every $H\subseteq G(n,r)$ satisfying $\deg_H(v) \geq(1/2+\eps)\deg_{G(n,r)}(v)$ for every $v \in V(H)$.
Note that our degree sequence matches the conjectured one, except for the degrees of a vanishing proportion of the vertices (roughly~$2rn$ of them, those which happen to fall within distance $r$ of $\{0\}$ or $\{1\}$).
In fact, using \cref{thm:deg-seq}, the number of vertices which do not satisfy the desired condition can be reduced to roughly $rn$, and we allow for roughly $rn$ vertices to have an even lower degree than in \cite[Conjecture~1.3]{ELW24}.


\subsection{Open problems}

Motivated by \cref{thm:cycle}, which shows that preserving slightly more than half of the degree of \emph{most} vertices of $P_n^k$ guarantees a Hamilton cycle, we wonder whether this is still the case if we preserve slightly more than half of the degree of \emph{each} vertex, and we propose the following conjecture (which would imply the case $d=1$ of \cite[Conjecture~1.3]{ELW24}).

\begin{conjecture}\label{conj:main}
    Let\/ $n\geq4$ and\/ $k\geq2$ be integers.
    Let\/ $G\subseteq P_n^k$ with $\deg_G(v)\geq\deg_{P_n^k}(v)/2+2$.
    Then,\/ $G$ is Hamiltonian.
\end{conjecture}



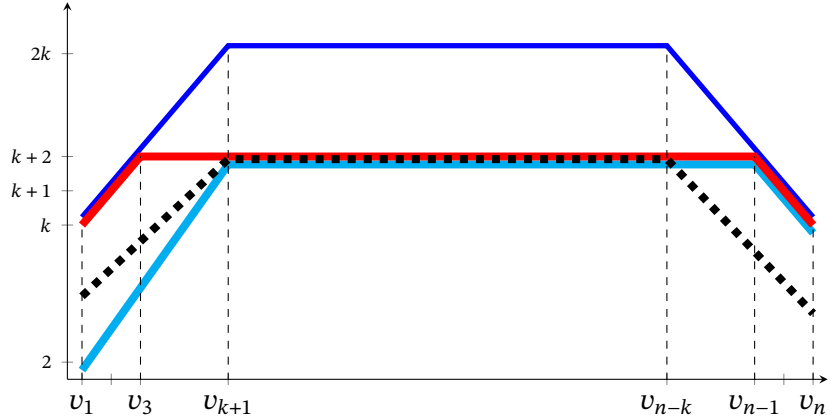
\begin{figure}
    \centering
    \pgfdeclarelayer{front}
    \begin{tikzpicture}
    \begin{axis}[
        width=10cm,
        height=5cm,
        scale only axis,
        axis lines=left,
        xmin=-0.5,
        xmax=25.5,
        xtick={0, 1, 2, 5, 20, 23, 24, 25},
        xticklabels={$v_1$,, $v_3$, $v_{k+1}$, $v_{n-k}$, $v_{n-1}$, ,$v_n$},
        ymin=1,
        ymax=23,
        ytick={2, 10, 12, 14, 20},
        yticklabels={\tiny{$2$}, \tiny{$k$}, \tiny{$k+1$}, \tiny{$k+2$}, \tiny{$2k$}},
    ]
        \addplot [yshift=3, blue, line width=2pt] coordinates { (0,10) (5,20) (20,20) (25,10) };
        
        \addplot [yshift=-3,cyan, line width=3pt
        ] coordinates { (0,2) (5,14)};
        
        \addplot [yshift=-3,cyan, line width=3pt
        ] coordinates {(5,14) (23,14) (24,12) (25,10) };     

        \addplot [yshift=0,red, line width=3pt] coordinates {(0,10) (1, 12) (2,14) (5,14)};  
        \addplot [yshift=0,red, line width=3pt] coordinates {(5,14) (23,14) (24,12) (25,10) };   

        \addplot [yshift=-1, black, dashed, line width=3pt
        ] coordinates {(0,6) (5,14) (20,14) (25,5)};

        \addplot [dashed, width=1pt] 
        coordinates {(0,0) (0,10)};
        \addplot [dashed, width=1pt] 
        coordinates {(2,0) (2,14)};
        \addplot [dashed, width=1pt] 
        coordinates {(25,0) (25,10)};
        \addplot [dashed, width=1pt] 
        coordinates {(23,0) (23,14)};
        \addplot [dashed, width=1pt] 
        coordinates {(5,0) (5,20)};
        \addplot [dashed, width=1pt] 
        coordinates {(20,0) (20,20)};

    \end{axis}
\end{tikzpicture}
    \caption{The blue line represents the degree of each vertex in $P_n^k$. The red (\cref{thm:cycle}) and cyan (\cref{thm:deg-seq}) lines indicate a sufficient minimum degree condition for Hamiltonicity.
    The dashed black line corresponds to \cref{conj:main}. }
    \label{fig:compare}
\end{figure}

In \Cref{fig:compare}, we show a comparison between the various degree conditions forcing the existence of a Hamilton cycle in \cref{thm:cycle,thm:deg-seq} as well as \cref{conj:main}.
This showcases the fact that \cref{thm:cycle} is weaker than \cref{conj:main}.
On the other hand, the result for Hamiltonicity from \cref{thm:deg-seq} and \cref{conj:main} are incomparable.

We remark that, if one could improve the degree sequence from \cref{thm:deg-seq} for $\ell=2$ so that the degree conditions on the last $k$ vertices mirror those of the first (which would be even stronger than \cref{conj:main}), then we would obtain \cref{conj:cycle} with an additive constant of $2$ instead of $1$ on the minimum degree condition.

In view of \cref{thm:cycle,thm:general}, it would be interesting to study the interplay between bandwidth and (effective) minimum degree as sufficient conditions for other spanning structures.
In particular, we suggest to investigate the containment of clique factors or of spanning bounded-degree trees.

\section{Proof of \texorpdfstring{\cref{thm:deg-seq}}{Theorem 4}}

    Given an oriented path, we call its first vertex its \emph{startpoint} and we call its last vertex its \emph{endpoint}.
    We allow paths to consist of a single vertex, in which case it is both the startpoint and the endpoint.
    Let $v_1,\ldots,v_n$ be a labelling of $V(G)$ which witnesses the fact that $\bw(G)\leq k$.
    Given any two distinct vertices $v_i,v_j\in V(G)$, we say that $v_i$ is \emph{to the left} of $v_j$ if $i<j$, and that it is \emph{to the right} of~$v_j$ otherwise.
    Let $i_1\coloneqq 1$, and let $i_2,\ldots,i_\ell$ denote the indices of the $\ell-1$ leftmost neighbours of $v_1$ (which exist since $\deg(v_1) \ge \ell$).
    Consider \cref{algo:general2}.

    \begin{algorithm}
    \caption{Construct a $K_{2,\ell}$-subdivision.}\label{algo:general2}
    \begin{algorithmic}[1]
    \State{For each $j\in[\ell]$, initialise an (oriented) path $P_{j}$ as the single vertex $v_{i_j}$.}
    \While{not all paths have crashed}
        \State{Let $\mathcal{P}\coloneqq\{P_j:j\in[\ell]\text{ and }P_j\text{ has not crashed yet}\}$.}
        \State{Let $v$ be the leftmost vertex among the endpoints of the paths in $\mathcal{P}$, and let $P$ be the path of~$\mathcal{P}$ whose endpoint is $v$.}
        \State{If $N(v)\subseteq \bigcup_{j\in[\ell]}V(P_j)$, say that $P$ has \emph{crashed} at $v$.}
        \State{Otherwise, let $v'$ be the leftmost vertex in $N(v)\setminus \bigcup_{j\in[\ell]}V(P_j)$ and extend $P$ by adding the (directed) edge $vv'$ to it.}
    \EndWhile
    \end{algorithmic}
    \end{algorithm}

    \Cref{algo:general2} produces $\ell$ vertex-disjoint paths. 
    We show now that they cover all the vertices of~$G$ and that their endpoints can be joined by edges in such a way that a subdivision of $K_{2,\ell}$ results.
    We begin by showing that each such path may only crash at one of the $\ell$ rightmost vertices.
    
    \begin{claim}\label{claim:crash22}
        None of the paths generated by \cref{algo:general2} may crash at a vertex $v_i$ with $i\leq n-\ell$.
    \end{claim}

    \begin{claimproof}

        First, note that the path starting at $v_1=v_{i_1}$ does not crash there, because $\deg_G(v_1) \ge \ell > \ell - 1$, which is the number of vertices initially covered by the other paths.
        Now consider any $2 \le i \le n-\ell$ and suppose for a contradiction that one of the paths, say~$P$, crashes at~$v_i$.
        We may assume that $v_i$ is the leftmost vertex where a path crashes.
        We analyse the situation at the instant when $P$ crashes at~$v_i$. Notice that, by the algorithm, $v_i$ is the leftmost endpoint at this time.
        Let $U\subseteq V(G)$ denote the set of right-neighbours of~$v_i$.
        Since $\bw(G)\leq k$ and $\deg(v_i) \geq \min\{i-1,k\}+\ell$, we know that $|U|\geq\ell$.
        As $P$ crashes at~$v_i$, the algorithm must already have covered all vertices in~$U$ with some path, and these vertices may have been used in different ways.
        In particular, there is at most one vertex $u'\in U$ such that $u'v_i\in E(P)$; if it exists, call it \emph{special} and observe that it must be the last vertex visited by~$P$ before reaching $v_i$.
        Moreover, since the algorithm maintains $\ell$ paths, there are at most $\ell-1$ vertices in $U$ which could be the endpoints of the other paths constructed so far.

        Let $U' \subseteq U$ be the (possibly empty) set of vertices of $U$ which are not the endpoints of a path nor the special vertex $u'$.
        Consider each $u \in U'$.
        Denote by $P_u$ the path in which $u$ is contained.
        Since $u$ is not the endpoint of $P_u$, the path $P_u$ must have been extended after reaching $u$.
        Since the algorithm did not append $v_i$ after $u$, the vertex following $u$ in $P_u$ must be some $\hat{u}$ to the left of $v_i$.
        In particular, since by assumption no path crashes at a vertex to the left of~$v_i$ and $v_i$ is the leftmost endpoint of a path, one of the following must occur:
        \begin{enumerate}[label=(\roman*)]
            \item\label{case1} at some point, $P_u$ jumps over $v_i$ again after having visited $\hat{u}$, or
            \item\label{case2} $P_u=P$ and, after having visited $\hat{u}$, $P$ reaches $v_i$ without ever jumping over it again.
        \end{enumerate}
        
        For each $u \in U'$ for which case~\ref{case1} holds, consider the first edge of $P_u$ that jumps over $v_i$ after having visited $\hat{u}$.
        This results in a collection $\mathcal{C}$ of vertex-disjoint edges of the form $wz$, with $w$ to the left of $v_i$ and $z$ to the right of $v_i$.
        If case~\ref{case2} never holds, then $|\mathcal{C}| = |U'| \ge |U|-\ell$.
        Note, moreover, that case~\ref{case2} can hold for at most one $u\in U'$ and, if it holds for any, then there is no special vertex.
        Thus, if case~\ref{case2} holds for some vertex, we have that $|\mathcal{C}| = |U'|-1 \ge |U| - (\ell-1)-1 = |U| - \ell$.
        In conclusion, it is always the case that $|\mathcal{C}| \ge |U| - \ell$.
        
        Additionally, observe that, the first time our algorithm extends any given path to a vertex to the right of~$v_i$ (which must have happened for at least one path, since $i>1$ and $v_i$ is the leftmost vertex where a path crashes), it does so with an edge $e$ disjoint from those in $\mathcal{C}$, as the edges of $\mathcal{C}$ belong to paths which visit a vertex of $U$ before such edge.
        Considering these edges in addition to those of $\mathcal{C}$ results in a collection $\mathcal{C}'$ of $|\mathcal{C}'|\geq|\mathcal{C}|+1\geq|U|-\ell+1$ vertex-disjoint edges with one vertex to the left of $v_i$ and the other to the right of $v_i$.


        Now let $W$ denote the set of endpoints to the left of $v_i$ of the edges in $\mathcal{C}'$.
        We must have $w\notin N(v_i)$ for all $w\in W$, as otherwise the algorithm would have chosen $v_i$ to extend some path from $w$ at an earlier step of the algorithm, rather than jumping over it.
        Moreover, since $\bw(G)\leq k$, each $w\in W$ and each neighbour of $v_i$ to its left correspond to some $v_j$ with $\max\{i-k,0\}<j<i$.
        As $\deg(v_i)-|U|$ is the number of left-neighbours of $v_i$ we then get
        \[
            \min\{i-1,k\} \ge |W|+\deg(v_i)-|U|\ge |W|+\min\{i-1,k\}+\ell-|U|,
        \]
        which implies that $|W|+\ell \le |U|$.
        This is a contradiction on the fact that $|W|=|\mathcal{C}'|\geq|U|-\ell+1$.
    \end{claimproof}

    It follows immediately from \cref{claim:crash22} that the endpoints of the paths constructed by \cref{algo:general2} are the vertices $v_i$ with $n-\ell<i\leq n$.
    Now we wish to prove that these paths contain all vertices of $G$.
    We say that a vertex $v$ is a \emph{gap} if it is not contained in any of the paths $P_i$ produced by \cref{algo:general2}.

    \begin{claim}\label{claim:gaps22}
        If $v_i$ is a gap, then there is a neighbour $v_j$ of $v_i$ to its right which is also a gap.
    \end{claim}

    \begin{claimproof}
        We argue similarly as above.
        By \cref{claim:crash22}, the vertices $v_i$ with $i>n-\ell$ are endpoints of paths, so they cannot be gaps.
        Now let $i\in[n-\ell]$, suppose that $v_i$ is a gap, let~$U$ denote the set of its right-neighbours, and suppose for a contradiction that none of them is a gap.
        Notice that, since~$v_i$ is a gap, none of the paths may have crashed at a vertex in $U$.
        Therefore, for every $u\in U$, when reached by a path, the algorithm must have chosen a neighbour to the left of $v_i$ to append to the path.
        In turn, this means that the union of the paths must contain a matching of size at least $|U|$ whose edges are of the form $wz$ with $w$ to the left of $v_i$ and $z$ to its right.
        Let $W$ denote the set of endpoints of these matching edges to the left of $v_i$, and observe that all $w \in W$ correspond to some $v_j$ with $\max\{i-k,0\}<j<i$. 
        Moreover, each $w\in W$ must not be a neighbour of~$v_i$, as otherwise $v_i$ would have been chosen by the algorithm when extending some path from $w$.
        This means that the number of neighbours of $v_i$ to its left is at most $\min\{i-1,k\}-|W|$.
        But~$v_i$ has at least $\min\{i-1,k\}+\ell-|U|$ neighbours to its left, which leads to a contradiction since $|W|\geq|U|$.
    \end{claimproof}
    
    Combined, \cref{claim:crash22,claim:gaps22} ensure that \cref{algo:general2} results in a set of $\ell$ pairwise vertex-disjoint paths $P_1, \dots, P_{\ell}$ such that, for each $j\in[\ell]$, the startpoint of $P_j$ is~$v_{i_j}$, the endpoints are $v_{n-\ell+1},\ldots,v_n$ (in some order), and the paths together cover all the vertices of~$G$.
    (Indeed, \cref{claim:crash22} gives the existence of paths with the desired endpoints and, since none of the last $\ell$ vertices is a gap and every gap must have another gap to its right by \cref{claim:gaps22}, there cannot be any gaps at all.)
    As $v_1$ is joined by an edge to all $v_{i_j}$ with $j\in[\ell]\setminus\{1\}$ and $v_n$ is joined by an edge to all $v_i$ with $n-\ell<i<n$, this immediately yields a subdivision of $K_{2,\ell}$ rooted at $v_1$ and $v_n$.\hfill\qedsymbol

\section*{Acknowledgement}

The research presented in this note was initiated during the Workshop \emph{Open Problems in Discrete Mathematics} organised by Anusch Taraz, Dennis Clemens, Fabian Hamann, Marco Wolkner, and Yannick Mogge. We would like to thank the organisers for their hospitality and the stimulating research environment that they created.

\bibliographystyle{mystyle} 
\bibliography{bibliography}

\end{document}